\documentclass[reqno,10pt]{amsart}
\usepackage{amssymb}
\usepackage{amsmath}
\usepackage{amscd}
\usepackage{amsthm, amsfonts,  amsmath, amssymb}
\usepackage{url}
\usepackage{amsmath, amssymb, amsthm, latexsym, enumerate}
\usepackage[top=2.9cm, bottom=2.5cm, left=3cm, right=3cm]{geometry}
\usepackage{geometry}
\usepackage{etoolbox}
\usepackage[all]{xy}
\DeclareFontEncoding{OT2}{}{} 

\newtheorem{prop}{Proposition}[section]
\newtheorem{defi}[prop]{Definition}

\newtheorem{conj}[prop]{Conjecture}

\newtheorem{lem}[prop]{ Lemma}
\newtheorem{thm}[prop]{Theorem}
\newtheorem{cor}[prop]{Corollary}

\newtheorem{remar}[prop]{Remark}

\newcommand{\QQ}{{\mathbb Q}}

\newcommand{\Q}{{\mathbb Q}}
\newcommand{\ZZ}{{\mathbb Z}}
\newcommand{\NN}{{\mathbb N}}

\newcommand{\OOO}{{\mathcal O}}

\newcommand{\MM}{{\mathfrak M}}

\newcommand{\BB}{{\mathfrak B}}

 \DeclareFontFamily{U}{wncy}{}
    \DeclareFontShape{U}{wncy}{m}{n}{<->wncyr10}{}
    \DeclareSymbolFont{mcy}{U}{wncy}{m}{n}
    \DeclareMathSymbol{\Sh}{\mathord}{mcy}{"58} 

\begin{document}
\title[A note on the Fourier coefficients of a Cohen-Eisenstein series]{\small {A note on the Fourier coefficients of a Cohen-Eisenstein series}
}

\author{SRILAKSHMI KRISHNAMOORTHY}

\date{\today}
\subjclass{}
\address{Indian Institute of Technology Madras\\ Tamil Nadu, India.}
\email{srilakshmi@iitm.ac.in, \ lakshmi@mpim-bonn.mpg.de}
\begin{abstract}
We prove a formula for the coefficients of a weight $3/2$ Cohen-Eisenstein series of square-free level $N$.
This formula generalizes a result of Gross and in particular, it proves a conjecture of Quattrini.
 Let $l$ be an odd prime number. For any elliptic curve $E$ defined over $\Q$   of rank zero and square-free conductor $N$, if $l \mid |E(\Q)| $, under certain conditions on the Shafarevich-Tate group $\Sh_D$, 
we show that $l$ divides  $|\Sh_D|$ if and only if  $l$ divides the class number $h(-D)$ of
$\Q(\sqrt{-D}).$
\end{abstract}

\maketitle
Keywords: Half-integral weight modular forms; Fourier coefficients, Shafarevich-Tate group.\\

Mathematics Subject Classification: Primary 11F37; Secondary: 11F67, 11R52.
\section{{  Introduction }}
Let $E$ be an elliptic curve of prime conductor $N$ and  analytic rank 0. Let $f$ be the new form of weight 2 of level $N$ on $\Gamma_0(N)$ associated to $E.$ Gross (Section 12, \cite{Gr87}) constructed 
$$\mathcal{G}=\sum_{D} m_D q^D,$$ a weight 3/2 modular form interms of certain modular forms $g_i$,  associated with $f.$
A special case of  Waldspurger's formula (Proposition 13.5,
\cite{Gr87}) relates the product of the $L$-functions $$L(f,1)L(f \times \epsilon_{-D},1)$$ to $m^2_D,$ where $(\frac{-D}{N}) \mathrm{sgn}(W_N) \neq -1$ and $f \times \epsilon_{-D}$ is the
cusp form corresponding to the twist by $-D$ of $E,$ and
$W_N$ is the Atkin-Lehner involution. Bocherer and Schulze-Pillot generalized
Gross's construction (Section 3, \cite{BS90}) for square-free level $N.$
Quattrini collected many numerical examples (Section 3.7, \cite{Qu11}) of certain definite quaternion algebras ramified at exactly one prime and presented a conjecture (Conjecture \ref{main}) on the coefficients of the Cohen-Eisenstein series $$\mathcal{H} = \sum^{n}_{i=1} \frac{1}{w_i} g_i,$$ where
$g_i$ and $w_i$ are certain orders defined in Section \ref{prelims}.
We work with certain definite quaternion algebras ramified at finitely many primes $p_1,$ $p_2,$ ...$p_k,$ and we compute the coefficients of $\mathcal{H}$ for square-free level
in theorem  \ref{mainthm}. As a consequence, we deduce the conjecture \ref{main} (Corollary \ref{main-conj}).

The estimation of the number of imaginary quadratic fields whose ideal class group has an element of order $l \geq 2$ and
the analogous questions for quadratic twists of elliptic curves has been the center of interest in many results. 
For elliptic curves $E$ of prime conductors, using the theory of $p$-adic $L$-functions and Eisenstein quotients, 
Mazur \cite{Ma79} showed that under certain conditions, the quadratic twist of $E$ by a primitive, odd quadratic 
Dirichlet character $\chi$ has finite Mordell-Weil group of order not divisible by a prime $l$ if and only if 
the quadratic field associated to $\chi$ has class number prime to $l.$
In \cite{Fr88}, Frey obtained the information about the elements of order $l$ in the Selmer group of $E_{D}$, the quadratic twist of $E$ by $-D$, by assuming
the elliptic curve $E$ over $\Q$ contains a $\Q$-rational torsion point of prime order $l$.
In \cite{Ja99}, James proved that  3 divides the order of the Selmer group of ${X_0(11)}_{D}$ if and only if 3 divides the class number $h(-D)$ under the similar assumption that
the elliptic curve $E$ contains a rational torsion point of order 3.
In \cite{Wo99}, Wong showed that there are infinitely many negative fundamental discriminants $-D$ such that the twist
$X_0(11)_{D}$ of the modular curve $X_0(11)$  has rank 0 over $\Q$ and an element of order 5 in its Shafarevich-Tate group.
Using the circle method and results of Frey, Kolyvagin,  Ono \cite{Ono01} proved a result for the nontriviality of class groups of imaginary quadratic fields and results on the nontriviality
of the Shafarevich-Tate groups of certain elliptic curves. It is also known that for almost all primes $l$, there exist infinitely many
square-free integers $D$ such that $l$ $\nmid$ $|\Sh_D|$ (\cite{Ko99}).

We prove that (Theorem \ref{keyprop}) if $E$ is an elliptic curve with square-free conductor $N$ and $l$ is an odd prime dividing $|E(\Q)| $, under certain conditions on the Shafarevich-Tate group $\Sh_D$,
the proportion of $\Sh_D$ in the family, divisible by $l$, is the same as the proportion of class numbers $h(-D)$ divisible by $l$ in the family of
negative quadratic fields $\Q(\sqrt{-D})$ with the same Kronecker conditions.

To prove theorem \ref{mainthm}, we follow  the strategy of Gross and we use Eichler's formula. The contents of this paper are as follows. In section \ref{prelims}, we discuss some preliminaries. In section \ref{optimal}, we compute the Fourier coefficients
of the modular forms $g_i$ in terms of $h(\mathcal{O}_{-D},R_i),$ the number of all optimal embeddings of the order of discrimiant $D$ into certain maximal orders $R_i$.
In section \ref{The order of the Shafarevich-Tate group}, we show that a certain odd prime divides the order of Shafarevich-Tate group of quadratic twists of elliptic curves  
if and only if it divides the class number of the corresponding imaginary quadratic field.
In section \ref{Cohen-Eisenstein series}, we compute the coefficients of the Cohen-Eisenstein series and we deduce Conjecture \ref{main}. 
\section{ {Preliminaries and statement of results}}\label{prelims}
Bocherer and Schulze-Pillot generalized
Gross's construction (Section 3, \cite{BS90}) for square-free level $N$ as follows. Let $B$ be a definite quaternion algebra ramified at primes $p_1, p_2,...,p_k$ and at $\infty$.
Let $N=p_1p_2...p_kM$ $(p_i \nmid M)$ be a square-free integer.
Let $\OOO$ be an order of level $N$.
Let $I_1$,$I_2$,...,$I_n$ be a set of left ideals representing  the distinct ideal classes of $\OOO$, with $I_1 = \OOO$.
Let $R_1$,$R_2$,...,$R_n$ be the respective right orders (of level $N$) of each ideal $I_i$.
For each $R_i$, let $L_i$ be the rank 3 lattice $\ZZ +2 R_i$. Denote the trace zero elements of $L_i$ by $S^{0}_i.$
For $b$ $\in$  $S^0_i,$ let $\NN(b)$ be the norm of $b.$ Let $w_i$ be the order of the finite group  $R_i^{*}/\pm 1$ for $i = 1$ to $n$. 
Define $$g_i = \frac{1}{2} \sum_{b \in S^0_i} q^{\NN(b)}.$$ The forms $g_i$ are in the Kohnen plus-space which is the space of modular forms $\sum a_n q^n$ of weight $3/2$ on  $\Gamma_0(4N)$
whose Fourier coefficients $a_n$ are 0 if  $-n \equiv 2, 3 \pmod 4.$

\subsection{Brandt matrices and Theta series}
Let $m$ be a positive integer. The Brandt matrix $B_m$ is defined by $B_m = (b_{ij}(m))_{n \times n},$ where $b_{ij}(m) = \frac{1}{e_j}  |\{ \alpha \in I^{-1}_j I_i \ : \ N(\alpha) \frac{N(I_j)}{N(I_i)} = m\}|,$
where $e_j = | R^{*}_j |.$\\
\hspace{0.25cm} The sum of any row in the matrix $B_m$ is given by $$b_m = \sum^{n}_{j=1} b_{ij}(m) = \sum_{d \mid m, (d, \frac{N}{M})=1 } d.$$
It is also the $m$-th coefficient of the zeta function $$\zeta_{\mathcal{O}} = \sum_{I} \frac{1}{{\mathbb{N}(I)}^{2s}} = \sum^{\infty}_{n=1}
\frac{b_n}{n^{2s}},$$ where the sum runs over all integral $\mathcal{O}$-left ideals $I.$
The vector $u=(1,1,...,1)$ is an eigenvector of the Brandt matrices, we have $B_m u^t = b_m u^t,$ for all positive integers $m.$
Fix $1 \leq i,j \leq n.$

These Bradnt matrices define a collection of theta series $$\theta_{ij}(\tau) = \frac{1}{e_j} \sum_{x \in I^{-1}_j I_i } q^{\frac{\NN(x) \NN(I_j)}{\NN(I_i)} }= \sum^{\infty}_{m=0} b_{ij}(m) q^m$$
which are modular forms of weight 2 and level $N.$

The series $e_2(z) = \sum^{n}_{i=1} \frac{1}{2w_i} + \sum^{\infty}_{m=1} b_m q^m$ is an Eisenstein series of weight 2 and level $N.$\\
Let $f$ be a new form of square-free level  $N=PM$ with $P=p_1p_2...p_k$ on $\Gamma_0(N)$  such that 
\begin{equation}\label{sgn-eqn1}
k \ \mathrm{is} \ \mathrm{odd,} \ 
\mathrm{sgn}(W_p) = -1, \ \mathrm{if} \ p \ \mid P, \ \mathrm{ sgn}(W_q) = +1, \ \mathrm{if} \ q \ \mid M.
\end{equation}
Suppose the elliptic curve $E$ corresponding to $f$ has analytic rank 0 and $l$ is an odd prime dividing the order of the torsion group of $E,$ then it can be shown that (Proposition 3.2, \cite{Qu11})
\begin{equation}\label{eqn10}
f \equiv e_2 \pmod l.
\end{equation}
\subsection{Waldspurger's formula}

The Shimura correspondence \cite{Sh73} relates the modular forms of half integral weight $k+1/2$ with classical modular forms of even weight $2k.$ We will define the modular form
$\mathcal{G}$ of weight $3/2$ which corresponds to the new form $f$ satisfying (\ref{sgn-eqn1}).  
Consider the quaternion algebra $B$ ramified exactly at $\infty$ and at the primes $p_i$ $\mid$ $N,$ 
where $\mathrm{sgn}(W_{p_i}) = -1.$ The Brandt matrices $B_m$ act on the the vector space
$V$ of formal linear combinations $\sum^n_{i=1} c_i I_i,$ $c_i \in \mathbb{C}.$ By Eichler's trace formula  there is a one to
one correspondence between Hecke eigenforms of weight 2 and level $N$
and eigenvectors in
$V$
of all Brandt matrices (up to a constant multiple) (Section 2, \cite{Po09}). Hence the normalized new form $f \in S_2(N)$ corresponds to a one-dimensional eigenspace $\langle v= ( v_1, v_2,..., v_n ) \rangle,$
of the Brandt matrices $\{{B}_p\}$ (of level $N$ and prime degree $p$) in $B$, such that ${B}_p v^{t} = a_p v^t,$
where $a_p$ is the eigenvalue satisfying $T_p f = a_p f,$ for all $p.$ We can assume that $(\frac{v_1}{w_1}, \frac{v_2}{w_2},...,\frac{v_n}{w_n})$ is primitive and has integer coordinates. Then
$$\mathcal{G} = \sum^{n}_{i=1} \frac{v_i}{w_i} g_i = \sum_{D}  m_D q^ D$$ is the weight $3/2$
modular form which corresponds to $f$ via the Shimura correspondence.\\
Let $P = p_1p_2...p_k.$ The modular form $\mathcal{G}$ is zero unless
$ \mathrm{sgn}(W_p) = \begin{cases}
  -1,  & \mbox{for } p \mid P \\
  +1, & \mbox{for } p \mid M \\
 \end{cases}$

If $-D$ is a fundamental discriminant such that $(\frac{-D}{p}) \mathrm{sgn}(W_p) \neq -1$ for every prime $p$
$\mid$ $\frac{N}{\mathrm{gcd}(N,D)},$ then the following
special case of Waldspurger's formula \cite{Wa81} holds (Section 3, \cite{BS90}).
\begin{equation}\label{eqn0}
\prod_{ p\mid  \frac{N}{\mathrm{gcd}(N,D)} }( 1+ (\frac{-D}{p}) sgn (W_p) ) L(f,1)L(f \otimes \epsilon_{-D}, 1) =
\frac{ 2^{\omega(N)} (f, f) m^2_D } {  \sqrt{D} \sum  \frac{ v^2_i }{ w^2_i}  }.
\end{equation}
\begin{defi}
The Cohen-Eisenstein series is the Eisenstein series of weight 3/2 corresponding to the eigenvector
$u=( 1, 1,..., 1 ),$ $\mathcal{H} := \sum^{n}_{i=1} \frac{1}{w_i} g_i.$
\end{defi}

\begin{remar}$(\mathrm{Multiplicity} \ \mathrm{one} \ \mathrm{modulo} \ l)$.\label{mul-remar}
When $N$ is prime, using the results of Mazur and Emerton \cite{Ma77}, \cite{Em02}, one can show that  the Brandt matrices $\{B_p \}$ reduced modulo $l$  have a dimension one eigenspace for the eigenvalues $\sigma(p)_N$
(Theorem 3.6, \cite{Qu11}).
Since $u$ and $v$ are both eigenvectors for the Brandt matrices $\{ B_p \}$, we have $\lambda u \equiv v \pmod l$
for some $\lambda \in \mathbb{F}^{\times}_l$.
If $N$ is square-free, then it is not clear whether the eigenspace corresponding to $u = (1,1,...,1)$  is one dimensional modulo $l$.
\end{remar}

Let $D$ be a natural number and let $\OOO_{-D}$ be the ring of integers in $\Q(\sqrt{-D})$. Let $h(-D)$ be the cardinality of the group $\mathrm{Pic}(\OOO_{-D})$, and
let $2u(-D)$ be the cardinality of the unit group $\OOO^{*}_{-D}.$
When $N$ is a square-free number, Quattrini made the following conjecture by observations on known congruences among weight two modular forms and known congruences
among eigenvectors of Brandt matrices. The details can be found in Section 3.1 -- 3.5 of \cite{Qu11}. 
\begin{conj}$(\mathrm{Conjecture}\ 3.7, \ $\cite{Qu11}$)$\label{main}Let $B$ be a definite quaternion algebra ramified at exactly one finite prime $p$ and let $N=pM$ $( \ p \nmid M \ )$
be a square-free integer. Let $ \mathcal{H} = \sum^{n}_{i=1} \frac{1}{w_i} g_i =  \sum^{n}_{i =1} \frac{ 1}{2w_i} +  \sum_{D > 0} \mathcal{H}(D) q^D.$
Let $D \in \NN$ be such that $-D$ is a fundamental discriminant and $( \frac{-D}{p}  ) \neq 1$, and
$( \frac{-D}{q}  )  \neq -1$ for every prime $q$ $\mid$ $M$.
Then
$$\mathcal{H}(D) = \frac{2^{\omega(N)-1 - s(D)} h(-D)}{u(-D)},$$
\end{conj}
where $\omega(N)$ is the number of distinct primes that divide $N$ and $s(D)$
is the number of primes that divide $N$ and ramify in $\QQ(\sqrt{-D}).$ 
If $M=1,$ then the above conjecture is the following result of Gross (Section 1, \cite{Gr87}).
\begin{prop}\label{gross-prop}
If $B$ is a definite quaternion algebra ramified only at a prime $N$ and $-D$ is a fundamental discriminant such that $( \frac{-D}{N}  ) \neq 1,$ then the coefficients  $\mathcal{H}(D)$ of the weight $3/2$ Eisenstein series
are given by  $$\mathcal{H}(D) = \frac{( 1 - (\frac{-D}{N}) )}{2}\frac{ h(-D)}{u(-D)}.$$
\end{prop}

In Conjecture \ref{main} and in Proposition  \ref{gross-prop}, Gross and Quattrini considered definite quaternion algebras ramified
at exactly one prime $p$ and at $\infty.$ We consider the generalized case, square-free level and definite quaternion algebras
ramified at finitely many primes $p_1,$ $p_2$,...$p_k$ and at $\infty$ (See Theorem \ref{mainthm}).
From Cremona's tables, the strong Weil curves of rank zero and
prime conductor with an odd torsion point, are listed by $E =
11A1$, $E=19A1$ and $E=37B1$. The first one has a 5-torsion point.
The other two curves have a 3-torsion point. For the $(-D)$
quadratic twists of $E$,  $ | \Sh_{D} |$ is  ${m^2_D}$, up to a
power of 2 and we also have $\lambda u \equiv v \pmod l,$ for some $\lambda \in \mathbb{F}^{\times}_l$
(remark \ref{mul-remar}). We state the following result of Quattrini (Proposition 3.8, \cite{Qu11}).
\begin{prop}
 Let $E$ be the strong Weil curve of rank 0 and prime conductor $N$.
 Consider the family $\{ E_D \}$ of negative quadratic twists of $E,$ for 
 $-D$ a fundamental discriminant and satisfying $(\frac{-D}{N}) =1.$  Suppose $E$ has a torsion point defined over $\Q,$ of odd prime order $l.$ Then, $|\Sh_D|$ is divisible by $l,$ if and only if the class number $h(-D)$ of $\Q(\sqrt{-D})$ is divisible by
 $l.$
 \end{prop}

We generalize the above proposition to square-free level $N$ as follows.

Let $E$ be an elliptic curve of analytic rank zero and square-free conductor $N=PM$ with $P=p_1p_2...p_k.$ 
Let $f$ be the new form of level $N$ on $\Gamma_0(N)$ corresponding to $E$ satisfying 
\begin{equation}\label{sgn-eqn}
\mathrm{Assume} \ k \ \mathrm{is} \ \mathrm{odd,} \ 
\mathrm{sgn}(W_p) = -1, \ \mathrm{if} \ p \ \mid P, \ \mathrm{ sgn}(W_q) = +1, \ \mathrm{if} \ q \ \mid M.
\end{equation}
 Consider the  family $\{ E_D \}$ of negative quadratic twists of $E$ satisfying the Kronecker condition
\begin{equation}\label{kronecker-eqn}
 ( \frac{-D}{p}  ) \neq 1 \ \mathrm{ for} \ p \ \mid P, ( \frac{-D}{q}  )  \neq -1, \ \mathrm{ for} \ q \ \mid M.
\end{equation}
We consider the definite quaternion algebra $B$ ramified exactly at all $p$ $\mid$ $P$ and at $\infty.$
We assume the following.
\begin{equation}\label{w-eqn}
 \mathrm{If} \ P \ \mathrm{is} \ \mathrm{composite}, \  \mathrm{then}  \ w_i \in \mathbb{F}^{\times}_l \ \mathrm{for} \ l = 3, 5 \ \mathrm{or} \ 7.
\end{equation}
The new form $f$ and the Eisenstein series $e_2$ of weight 2 correspond to the 3/2 weight forms 
$\mathcal{G}$ and the Cohen-Eisenstein series $\mathcal{H}$ respectively, under the Shimura correspondence.
Let $v$ and $u$ be the eigenvectors of the Brandt matrices associated with the forms $f$ and $e_2$  respectively.
Suppose $\lambda u \equiv v \pmod l$
for some $\lambda \in \mathbb{F}^{\times}_l,$ then the congruence (\ref{eqn10}) in weight 2
can be lifted to a congruence  in weight 3/2,
\begin{equation}\label{eqn11}
 \lambda \mathcal{G} \equiv \mathcal{H} \pmod l.
\end{equation}
Thus we have the following result

\begin{thm}\label{keyprop}
 Let $E$ be an elliptic curve of analytic rank zero and square-free conductor $N=PM.$ Let $f$ be the new form of level
 $N$ corresponding to $E$ satisfying (\ref{sgn-eqn}).  
 Consider the  family $\{ E_D \}$ of negative quadratic twists of $E$ satisfying the Kronecker condition
 (\ref{kronecker-eqn}). Suppose $E$ has a torsion point defined over $\Q,$ of odd prime order $l$ and that $|\Sh_{D}| =
 {m^2_D}$ (upto a power of 2).  Assume that $\lambda u \equiv v \pmod l,$ for some $\lambda \in \mathbb{F}^{\times}_l$ and
 (\ref{w-eqn}) holds. Then, $|\Sh_D|$ is divisible by $l,$ if and only if the class number $h(-D)$ of $\Q(\sqrt{-D})$ is divisible by
 $l.$
\end{thm}


\section{Optimal embeddings}\label{optimal}
We continue with the notation set out in the previous sections.
Let $K$ be a quadratic field over $\Q.$
Let $\phi$ be an embedding of $K$ into $B$.
The field $K$ is totally imaginary as $B$ is a definite quaternion algebra. Let $\OOO_{-D}$ be an order of $K$ of discriminant $D$.
\begin{defi}
We say that $\phi$ is an  optimal embedding of the order $\OOO_{-D}$ into $R_i$ if $\phi$ is an embedding of $K$ into $B$ such that
$\phi(\OOO_{-D}) = \phi(K) \cap R_i$.
\end{defi}
Two optimal embeddings $i_1, i_2$ are equivalent if they are conjugate to each other by
an element in $R^{*}_i$. In other words, if there exists $x$ $\in$ $R^{*}_i$ such that $i_1(y) = x i_2(y) x^{-1}$
for all $y$ $\in$ $K$.

The Legendre symbol $( \frac{-D}{p} )$ is defined by
$( \frac{-D}{p} ):=$
$\begin{cases}
  1,  & \mbox{if } p \mbox{ splits in } K \\
  0, & \mbox{if }  p  \mbox{ ramifies in } K \\
 -1, & \mbox{if }  p  \mbox{ is inert in } K.
\end{cases}$

The Eichler symbol $\{ \frac{-D}{p} \}$ is defined by
 $\{ \frac{-D}{p} \}:=$
 $\begin{cases}
  1,  &  \mbox{if} \  {p^2}  \mid  D \\
  0,  &  \mbox{if}  \ {p} \mid D, {p^2} \nmid {D}\\
  ( \frac{-D}{p} ), & \mbox{if } p \nmid D .
\end{cases}$
\\
\\
We prove a lemma and a proposition. We will use them in the proof of Theorem \ref{mainthm}.
\begin{lem}\label{for1} Let $h(\OOO_{-D}, R_i)$ be the number of equivalence classes of
 optimal embeddings of the order of discriminant $D$ into $R_i$.
Then
$$ \sum^{n}_{i=1} h(\OOO_{-D}, R_i) =  h(-D) \prod^{k}_{i=1} ( 1-  \{ \frac{-D}{p_i} \} ) \prod_{q \mid M}  ( 1+  \{ \frac{-D}{q} \} ).$$
\end{lem}
\begin{proof}
Let $\{\MM\}$ be a system of representatives of two-sided $R_i$ ideals modulo two-sided $R_i$ ideals of the form $R_i \xi $ where $\xi$ is an $\OOO_{-D}$
ideal. Let $\{ \BB \}$ be a system of representatives of the ideal classes in $\OOO_{-D}.$
Consider the set of all $(\MM, \BB )$ such that\\
(1) The norm of $\MM$ is square-free and if $q$ is a prime divisor of the norm of $\MM$,
then either $q=p_i$ (for some $i =1$ to $k$) with $\{ \frac{-D}{p_i} \}=-1$ or $q$ is a prime divisor
of $M$ with $\{ \frac{-D}{q} \}=1$ and\\(2) $\BB$ is an integral ideal coprime to the conductor of $\OOO_{-D}$.\\It is easy to observe that the
number of $ (\MM,  \BB)$ satisfying (1) and (2) is equal to
$$h(-D) \prod^{k}_{i=1} ( 1-  \{ \frac{-D}{p_i} \} ) \prod_{q \mid M}  ( 1+  \{ \frac{-D}{q} \} ).$$
There is a one-to-one correspondence between the set of all $(\MM,  \BB)$ satisfying (1) and (2) and equivalence classes of optimal embeddings of the order of discriminant $-D$ into $R_i$.
For the proof of this correspondence, we refer to Section 3.2 of \cite{Sh65} (or) Satz 6,7 of \cite{Ei55}.
\end{proof}
We compute the Fourier coefficients of the modular forms $g_i,$ for $i=1$ to $k$ in the following proposition.
\begin{prop}\label{for2}
 Let $g_i =\frac{1}{2} +  \frac{1}{2}  \sum_{ D > 0} a_i(D) q^{D}.$ Then
 $a_i(D)$ is the number of elements $b$ $\in$ $R_i$ with $\mathrm{Tr}(b) =0$, $b \in \ZZ + 2 R_i$, $\NN(b) = D.$
 For $i = 1$ to $n$, we have $$a_i(D) = w_i \sum_{-D=df^2} \frac{h(\OOO_d, R_i)} {u(d)}, $$ where $u(d)=1$ unless $ d =-3, -4$ when $u(d) = 3,2$ respectively.
\end{prop}
\begin{proof}
Let $S$ be the set of elements  $b$ $\in$ $R_i$ with $\mathrm{Tr}(b) =0$, $b \in \ZZ + 2 R_i$ and $\NN(b) = D.$\\
For a negative integer $d$,
if $f: \QQ(\sqrt{d}) \hookrightarrow B$ is an embedding of an order $\OOO_d$ into $R_i$, then\\
$b=f(\sqrt{d})$ is an element with trace 0 and norm $-d$. Since $\OOO_{d} = \ZZ + \ZZ \frac{(-d + \sqrt{d})}{2},$ we have $b$ $\in$ $(\ZZ + 2R_i)$.
Hence $b$ $\in$ $S^{0}_{i} = \{ x  \in B | \mathrm{Tr}(x) =0 \} \cap (\ZZ + 2R_i).$\\ Conversely, if $b$ is an element in $S^{0}_i$ with norm $-d$ , then  $f(\sqrt{d}) = b$ gives rise to an embedding
of the order $\OOO_{d} = \ZZ + \ZZ \frac{ (-d + \sqrt{d})}{2}$ into $R_i$. The embedding $f(\sqrt{d})=b$ is optimal if and only if $b \notin f(\ZZ + 2R_i)$ for some $f >1.$
Let $h^{*}(\OOO_{-D}, R_i)$ be the the number of optimal embeddings of $\OOO_{-D}$ into $R_i$.
Using the above connection we proved that
$$a_{i}(D) = |S| = \sum_{-D=df^2} \{ b \in S,  \frac{b}{f} \in S^{0}_{i}, \frac{b}{f} \notin n(\ZZ+2R_i)  \ \mathrm{for} \ n>1 \} = \sum_{-D=df^2} h^{*}(\OOO_d, R_i).$$
The group $\Gamma_i = R^{*}_i/\pm1$ acts on $S$. The $\Gamma_i$ orbits of $S$ correspond to equivalence classes of optimal embeddings.
Hence $$ | S/\Gamma_i | =  \sum_{-D=df^2} h(\OOO_d, R_i).$$ The order of the stabilizer of an element $b \in S$ is 1 unless
the corresponding embedding extends to $\ZZ[\mu_6]$ or $\ZZ[\mu_4]$, when it is 3 or 2 respectively. Thus we have shown that
$$a_i(D) = w_i \sum_{-D=df^2} \frac{h(\OOO_d, R_i)} {u(d)}, $$ where $w_i = |\Gamma_i|.$
\end{proof}
Gross computed the traces of the Brandt matrices for prime level case\\(cf. Proposition 1.9, \cite{Gr87}). 
It holds for square-free level, as we state in the following.
\begin{prop}
For all $m \geq 0,$
$$\mathrm{Tr}(B(m)) = \sum_{s \in \mathbb{Z}, s^2-4m \leq 0} \mathcal{H}(4m-s^2).$$
\end{prop}
\begin{proof} The diagonal entry of the brandt matrix $B(m)$ is  $b_{ii}(m) = \frac{1}{e_i}|\{b, b \in R_i, \NN(b) = m\} |.$\\
If $m =0,$ then $$\mathrm{Tr}(B(0)) =  \frac{1}{24}
\prod^{k}_{i=1} ( {p_i} -1  ) \prod_{q \mid M}   ( q+1 ) =
\sum^{n}_{i =1} \frac{ 1}{2w_i} = \mathcal{H}(0).$$ Let $A_i(s,m)$
be the set of elements  $b$ $\in$ $R_i$ with $\mathrm{Tr}(b) =s$
and $\NN(b) = m.$\\ This is a finite set. If $s^2-4m > 0$, then it
is an empty set.  Hence
$$\mathrm{Tr}(B(m)) =  \sum^n_{i=1} b_{ii}(m) = \sum^n_{i=1}
\sum_{s^2 \leq 4m} \frac{|A_i(s,m)|}{|R_i^{*}|} = \sum_{s^2 \leq
4m} ( \sum^n_{i=1} \frac{|A_i(s,m)|}{|R_i^{*}|} ).$$ If $s^2 =
4m$, then the inner sum $$\sum^n_{i=1}
\frac{|A_i(s,m)|}{|R_i^{*}|} = \sum^{n}_{i =1} \frac{ 1}{2w_i} =
\mathcal{H}(0).$$ Assume that $D =4m-s^2
>0.$
 As in the
proof of Proposition \ref{for2}, we can show that
$$\frac{|A_i(s,m)|}{|R_i^{*}|} = \sum_{-D=df^2}
\frac{1}{2}\frac{h(\OOO_d, R_i)}{u(d)}.$$ By Lemma \ref{for1} and
Theorem \ref{mainthm}, $$\sum^n_{i=1} \frac{|A_i(s,m)|}{|R_i^{*}|}
= \sum^{n}_{i=1} \sum_{-D=df^2} \frac{1}{2}\frac{h(\OOO_d,
R_i)}{u(d)} = \mathcal{H}(4m-s^2).$$
\end{proof}

\section{The order of the Shafarevich-Tate group}\label{The order of the Shafarevich-Tate group}
Recall that we have equation (\ref{eqn0}) which relates the L-function of $f$ with the coefficients $m^2_D,$
$$\prod_{ p\mid  \frac{N}{\mathrm{gcd}(N,D)} }( 1+ (\frac{-D}{p}) sgn (W_p)  ) L(f,1)L(f \otimes \epsilon_{-D}, 1) =
\frac{ 2^{\omega(N)} (f, f ) m^2_D } {  \sqrt{D} \sum  \frac{ v^2_i }{ w^2_i}  }.$$
If $E$ is the elliptic curve with conductor $N$ associated with $f \in S_2(\Gamma_0(N)),$
then we have $L(E,1)=L(f,1).$  Then the $L$-function $L(f \otimes \epsilon_{-D}, 1) = L(E_{D},1),$
where $E_{D}$ is the $-D$ quadratic twist of $E$ associated with $f \otimes \epsilon_{-D} \in S_2(\Gamma_0(ND^2)).$
Assume that the rank of $E$ is $0.$ The rank $0$ case of Birch and Swinnerton-Dyer Conjecture gives  $$\frac{L(f \otimes \epsilon_{-D}, 1)}{\Omega_D}= \frac{L(E_{D},1)}{\Omega_D}= \frac{|\Sh_{D}| \prod c_{p,D}}{|\mathrm{Tor}(E_{D})|^2},$$
where $c_{p,D}$'s  are the Tamagawa numbers and $\mathrm{Tor}(E_{D})$ is the torsion subgroup of $E_{D}(\Q),$
${\Omega_D}$ is the real period of $E_D.$  Let $$C(D) = \frac{\prod_{ p\mid  \frac{N}{\mathrm{gcd}(N,D)} }( 1+ (\frac{-D}{p}) sgn (W_p)  ) } {2^{\omega(N)}}\frac{\Omega_D \prod c_{p,D}  \sqrt{D} \frac{ v^2_i }{ w^2_i}L(f,1)}{(f, f ) |\mathrm{Tor}(E_{D})|^2 }.$$
Then $|\Sh_D| = \frac{m^2_D}{C(D)}.$
Math softwares can be used to compute the term $C(D).$ 
\subsection{Proof of Theorem \ref{keyprop}}
\begin{proof}We prove the theorem  when $P$ is prime. One can conclude the theorem similarly when $P$ is composite.
If $l$ is an odd prime  dividing the order of the group of torsion points of the elliptic curve $E$,
by Mazur's theorem, $l = 3, 5$ or $7.$ We know that $w_i \mid 12,$ the product $\prod^{n}_{i=1} w_i$ equals the exact
denominator of $\frac{N-1}{12}$ and $3$ divides the exact numerator of $\frac{N-1}{12}.$
Hence $w_i  \in \mathbb{F}^{\times}_l$ for $l = 3, 5$ or $7.$
From $\lambda \mathcal{H} - \mathcal{G} = \sum^{n}_{i=1} \frac{(\lambda - v_i)}{w_i} g_i,$
it follows that the congruence $\lambda u \equiv v \pmod l,$ for some $\lambda \in \mathbb{F}^{\times}_l$  gives a congruence
$\lambda \mathcal{H} \equiv \mathcal{G} \pmod l.$
This yields a congruence on the coefficients $\lambda \mathcal{H}(D) \equiv m^2_D \pmod l.$  From Corollary \ref{main-cor}, we see that
$l$ divides $\mathcal{H}(D)$ if and only if $l$ divides $h(-D).$ We also have $|\Sh_{D}| =
 {m^2_D}$ (up to a power of 2). Hence $|\Sh_D|$ is divisible by $l,$ if and only if the class number $h(-D)$ of $\Q(\sqrt{-D})$ is divisible by
 $l.$ 
\end{proof}

By letting $k=1$ in the above Theorem, we deduce the following corollary.

\begin{cor}[Proposition 3.9, \cite{Qu11}] Let $E,$ $E_D,$ $l$ and $\Sh_D$ be as in Theorem \ref{keyprop}.
Assume that there is exactly one prime $p \mid N$ such that the sign of $ W_p = -1.$
Then, $|\Sh_D|$ is divisible by $l,$ if and only if the class number $h(-D)$ of $\Q(\sqrt{-D})$ is divisible by  $l.$
\end{cor}
\section{Cohen-Eisenstein series}\label{Cohen-Eisenstein series}

\subsection{Examples}
We calculate the Fourier coeffients of the weight 3/2 Eisenstein series $\mathcal{H} = \sum^{n}_{i=1} \frac{1}{w_i} g_i$ 
and the class numbers of imaginary quadratic fields $K=\Q(\sqrt{-D})$ for $d \leq 2000$ by using MAGMA. 
Let $D > 0$ be a natural number and let $\OOO_{-D}$ be the ring of integers in $\Q(\sqrt{-D})$. Let $h(-D)$ be the cardinality of the group $\mathrm{Pic}(\OOO_{-D})$, and
let $2u(-D)$ be the cardinality of the unit group $\OOO^{*}_{-D}.$ A prime $l$ is inert, splits or ramifies in $\OOO^{*}_{-D}$
if the Kronecker symbol $( \frac{-D}{l} )$ is -1, 1, 0 respectively.

We consider the strong Weil curves of rank zero with an odd torsion point from Cremona's table \cite{Cr97}.
\\
$\bullet N=66 = 2.3.11$
\\
We have elliptic curve $E = 66 \mathrm{C(I)} = [1,0,0,-45,81]$ of level $66$ with  analytic rank zero and $|\mathrm{Tor}(E)| =10.$
We have $\mathrm{sgn} (W_2) = \mathrm{sgn} (W_3)  = \mathrm{sgn} (W_{11})  = -1.$
We work in the quaternion algebra ramified at 2, 3, 11 and at $\infty.$
We calculate the Brandt matrices for an order of level 66.  We have,
for $D \leq 2000$ such that $-D$ is a fundamental discriminant and $( \frac{-D}{2} ), ( \frac{-D}{3} )$ and $( \frac{-D}{11} ) \neq 1:$
\\
\\
$\mathcal{H}(D):=$
$\begin{cases}
 \frac{2^{2} h(-D)}{u(-D)} ,  & \mbox{if none of the primes} \  2,3, 11 \ \mbox{ramifies in} \  K  \\
  \frac{2^{1} h(-D)}{u(-D)}, & \mbox{if  exactly one prime} \ p \mid 66  \ \mbox{ramifies in  } K \\
 \frac{h(-D)}{u(-D)}, & \mbox{if  exactly two primes} \ p \mid 66 \ \mbox{ramify in  } K\\
 \frac{ h(-D)}{2 u(-D)}, & \mbox{if} \ 2, 3 \ \mbox{and} \ 11 \ \mbox{ramify in  } K.
\end{cases}$
\\
\\
$\bullet N=210 = 2.3.5.7$
\\
We have elliptic curve $E = 210 \mathrm{A(A)} = [1,0,0,-41,-39]$ of level $210$ with  analytic rank zero and $|\mathrm{Tor}(E)| = 6.$
We have $\mathrm{sgn} (W_2) = \mathrm{sgn} (W_3)  = \mathrm{sgn} (W_{7})  = -1$ and $\mathrm{sgn} (W_5) = +1.$ 
We work in the quaternion algebra ramified at 2, 3, 7 and at $\infty.$
We calculate the Brandt matrices for an order of level 210.  We have,
for $D \leq 2000$ such that $-D$ is a fundamental discriminant and $( \frac{-D}{2} ), ( \frac{-D}{3} ), ( \frac{-D}{7} ) \neq 1$ and
$( \frac{-D}{5} ) \neq -1:$
\\
\\
$\mathcal{H}(D):=$
$\begin{cases}
 \frac{2^{3} h(-D)}{u(-D)} ,  & \mbox{if none of the primes} \  2,3, 5, 7 \ \mbox{ramifies in} \ K \\
 \frac{2^{2} h(-D)}{u(-D)} ,  & \mbox{if exactly one of the primes} \  p \mid 210 \ \mbox{ramifies in} K  \\
  \frac{2^{1} h(-D)}{u(-D)}, & \mbox{if  exactly two of the primes} \ p \mid 210  \ \mbox{ramify in  } K \\
 \frac{h(-D)}{u(-D)}, & \mbox{if  exactly three primes} \ p \mid 210 \ \mbox{ramify in  } K\\
 \frac{ h(-D)}{2 u(-D)}, & \mbox{if} \ 2, 3, 5 \ \mbox{and} \ 7 \ \mbox{ramify in   } K.
\end{cases}$
\\
\\
We have also computed the Fourier coefficients $\mathcal{H}(D)$
for the rank 0 elliptic curves $E = 110A1(C) = [1,1,1,10,-45]$ with $\mathrm{Tor}(E) = 5,$
$E = 114A(A) = [1,0,0,-8,0]$ with  $\mathrm{Tor}(E) = 6,$
$E = 130B(A) = [1,-1,1,-7,-1]$ with  $\mathrm{Tor}(E) = 4,$
$E = 210B(A) = [1,0,1,-498,4228]$ with  $\mathrm{Tor}(E) = 6$
and several other examples. Based on our numerical examples, we observed a generalization 
of the conjecture \ref{main} which we prove in Corollary \ref{main-cor}.
\begin{thm}\label{mainthm}
Let $B$ be a definite quaternion algebra ramified at $p_1, p_2,...,p_k$. Let $N=p_1p_2...p_kM$ $( \ p_i \nmid M \ )$
be a square-free integer. Denote by $ \mathcal{H} = \sum^{n}_{i=1} \frac{1}{w_i} g_i =  \sum^{n}_{i =1} \frac{ 1}{2w_i} +
\sum_{D > 0} \mathcal{H}(D) q^D.$ Then we have
$$ \mathcal{H}(D) =\frac{1}{2}\sum_{ -D=df^2  }
 \Big[\frac{h(d)}{u(d)} \prod^{k}_{i=1} \Big( 1-  \{ \frac{d}{p_i} \} \Big) \prod_{q \mid M}  \Big( 1+  \{ \frac{d}{q} \} \Big) \Big].$$
\end{thm}
\begin{proof} Consider the weight $3/2$ Cohen-Eisenstein Series $\mathcal{H}$,
$$ \mathcal{H} = \sum^{n}_{i=1} \frac{1}{w_i} g_i = \frac{1}{2} \sum^{n}_{i=1} \frac{1}{w_i} +  \sum_{D >0}  \sum^{n}_{i=1} \frac{a_i(D)} {w_i} q^D.
$$
By Proposition \ref{for2}, we see that
\begin{equation}\label{eqn1}
\sum^{n}_{i=1} \sum_{D >0} \frac{a_i(D)}{w_i} q^D = \frac{1}{2} \sum_{ D >0} \sum_{ -D=df^2  }
\Big(\sum^{n}_{i=1} \frac{ h(\OOO_d, R_i)}{u(d)} \Big).
\end{equation}
By Lemma \ref{for1}, we have
\begin{equation}\label{eqn2}
\sum^{n}_{i=1} h(\OOO_d, R_i) =
h(d) \prod^{k}_{i=1} \Big( 1-  \{ \frac{d}{p_i} \} \Big) \prod_{q \mid M}  \Big( 1+  \{ \frac{d}{q} \} \Big).
\end{equation}
Substituting equations (\ref{eqn1}) and (\ref{eqn2}) in the Fourier expansion of $\mathcal{H},$ we get

$$\mathcal{H} = \sum^{n}_{i=1} \frac{1}{w_i} g_i = \frac{1}{2} \sum^{n}_{i=1} \frac{1}{w_i} + \frac{1}{2}
\sum_{ D>0} \sum_{ -D=df^2  }
 \Big[\frac{h(d)}{u(d)} \prod^{k}_{i=1} \Big( 1-  \{ \frac{d}{p_i} \} \Big) \prod_{q \mid M}  \Big( 1+  \{ \frac{d}{q} \} \Big) \Big] q^D.$$
Hence
\begin{equation}\label{main-eqn}
\mathcal{H}(D) = \frac{1}{2}\sum_{ -D=df^2  }
 \Big[\frac{h(d)}{u(d)} \prod^{k}_{i=1} \Big( 1-  \{ \frac{d}{p_i} \} \Big) \prod_{q \mid M}  \Big( 1+  \{ \frac{d}{q} \} \Big) \Big].
\end{equation}
 This completes the proof.
\end{proof}
 We deduce the following corollary which generalizes the result of Gross (Proposition \ref{gross-prop}) for square-free level $N.$ 
\begin{cor}\label{main-cor} Let $B$ and $\mathcal{H}$ be as in Theorem \ref{mainthm}.
If $-D$ is the fundamental discriminant, $\omega(N)$ is the number of distinct primes that divide $N$, $s(D)$
is the number of primes that divide $N$ and ramify in $\QQ(\sqrt{-D}),$ and
$( \frac{-D}{p_i}  ) \neq 1$, for every  $i = 1$ to $k$, and
$( \frac{-D}{q}  )  \neq -1$ for every prime $q$ $\mid$ $M$, then
$$\mathcal{H}(D) = \frac{2^{\omega(N)-1 - s(D)} h(-D)}{u(-D)}.$$
\end{cor}
\begin{proof} If $-D$ is the fundamental discriminant satisfying the Kronecker conditions, then from the equation (\ref{main-eqn}), we have
$$\mathcal{H}(D) =  \prod^{k}_{i=1} \Big( 1-  \{ \frac{-D}{p_i} \} \Big) \prod_{q \mid M}  \Big( 1+  \{ \frac{-D}{q} \} \Big) \frac{1}{2} 
 \frac{h(-D)}{u(-D)}= \frac{2^{\omega(N)-1 - s(D)} h(-D)}{u(-D)}.$$
\end{proof}

\begin{cor}\label{main-conj}
 Conjecture \ref{main} holds.
\end{cor}
\begin{proof}
 Conjecture \ref{main} follows immediately from Corollary \ref{main-cor}
by letting $k=1.$
\end{proof}

\section*{Acknowledgements}
This work was done when the author was a Visiting Postdoctoral Fellow at the Max Planck Institute for Mathematics, Bonn.
She would like to thank the institute for providing excellent working conditions. She would like to thank  Neil Dummigan, ‎Tomoyoshi Ibukiyama and Narasimha Kumar and for helpful discussions.
She would also like to thank  Chitrabhanu Chaudhuri, Rafael von K\"{a}nel, Rachel Newton and Sarang Sane for helpful comments  on earlier drafts of this paper.
The research of the author was supported by a DST-INSPIRE Grant.


\begin{thebibliography}{Bz}

\bibitem[BS90]{BS90} S. Bocherer, R. Schulze-Pillot, On a theorem of Waldspurger and on Eisenstein series of Klingen type,
{\em Math. Ann.}, 288, (1990), 361--388.
\bibitem[Ei55]{Ei55} M. Eichler, Zur Zahlentheorie der Quaternion-Algebren. {\em J. Reine Angew. Math.}, {\bf 195} (1955), 127--151.
\bibitem[Cr97]{Cr97} J.E. Cremona, Algorithms for Modular Elliptic Curves, Cambridge University Press, 1997.
\bibitem[Em02]{Em02}  M. Emerton, Supersingular elliptic curves, theta series and weight two modular forms.  {\em J. Amer. Math. Soc.} Vol {\bf 15}, (2002), no. 3, 671--714.
\bibitem[Fr88]{Fr88} G. Frey, On the Selmer group of twists of elliptic curves with Q-rational torsion points, {\em Canad. J. Math.}, Vol {\bf 40}, no. 3, (1988) 649--665.
\bibitem[Gr87]{Gr87} B. Gross, Heights and the special values of $L$-series, {\em CMS conference Proceedings}, Vol {\bf 7}, (1987),\\
115--187.
\bibitem[Ko99]{Ko99} W. Kohnen, K. Ono, Indivisibility of class numbers of imaginary quadratic fields and orders of Tate-Shafarevich groups of elliptic curves with complex multiplication.
{\em Invent. Math.} Vol {\bf 135} (1999), no. 2, 387--398. 
\bibitem[Ma77]{Ma77} B. Mazur,  Modular curves and the Eisenstein ideal. {\em Inst. Hautes Etudes Sci. Publ. Math.}, (1977), no. 47,  33--186.
\bibitem[Ma79]{Ma79} B. Mazur, On the arithmetic of special values of L functions. {\em Invent. Math.}, Vol {\bf 55} (1979), no. 3, 207--240. 
\bibitem[Po09]{Po09} P. Ponomarev, New forms of squarefree level and theta series. {\em Math. Ann.} Vol {\bf 345} (2009), no. 1, 185--193. 
\bibitem[Qu11]{Qu11} P.L. Quattrini, The effect of torsion on the distribution of $\Sh$ among quadratic twists of an elliptic curve, {\em Journal of Number Theory}, no. 2, (2011), 195--211.
\bibitem[Sh65]{Sh65} H. Shimizu, On Zeta Functions of Quaternion Algebras {\em Ann. Math.} Vol {\bf 81} no. 1, (1965), 166--193.
\bibitem[Sh73]{Sh73} G. Shimura, {\em On modular forms of half-integral weight}, Ann. of Math. {\bf 97} (1973), 440--481.
\bibitem[Ono01]{Ono01} K. Ono, Nonvanishing of quadratic twists of modular $l$-functions and applications to elliptic curves
{\em J. Reine Angew. Math.} Vol {\bf 553} (2001), 81--97.
\bibitem[Ja99]{Ja99} K. James, Elliptic curves satisfying the Birch and Swinnerton Dyer conjecture and mod 3 {\em J. Number Theory.}
 Vol {\bf 128} (2008) 2823--2835.
 \bibitem[Wa81]{Wa81} J.L. Waldspurger, Sur les coefficients de Fourier des formes modulaires de poids demi-entier. {\em J. Math. Pures Appl.} (9) Vol {\bf 60} (1981), no. 4, 375--484. 
 \bibitem[Wo99] {Wo99} S. Wong, Elliptic curves and class number divisibility {\em Int. Math. Res. Not.}, {\bf 12} (1999), 661--672.

\end{thebibliography}
\end{document}